\documentclass{article}

\newtheorem{theorem}{Theorem}[section]
\newtheorem{corollary}[theorem]{Corollary}

\newtheorem{lemma}[theorem]{Lemma}
\newtheorem{example}[theorem]{Example}

\usepackage{amsmath,amsfonts,amssymb}

\usepackage[all]{xy}
\input xy
\xyoption{all}

\begin{document}

\title
{Surjectivity of isometries of weighted 
spaces of holomorphic functions and of Bloch spaces}

\author{Christopher Boyd and Pilar Rueda\footnote{The second author was supported by Ministerio de 
Econom\'{\i}a yCompetitividad (Spain) MTM2011-22417.} }

\maketitle

\begin{abstract}
We examine the surjectivity of isometries between weighted spaces of holomorphic
functions. We show that for certain classical weights on the open unit disc all 
isometries of the weighted space of holomorphic functions, ${ \mathcal H}_{v_o}(
\Delta)$, are
surjective. Criteria for surjectivity of isometries of ${ \mathcal H}_v(U)$ in 
terms 
of a separation 
condition on points in the image of
${ \mathcal H}_{v_o}(U)$ are also 
given for $U$ a bounded open set in $\mathbb{C}$. Considering the weight $v(z)=
1-|z|^2$ and the isomorphism $f\mapsto f'$ we are able to show that all 
isometries of the little Bloch space are surjective.
\end{abstract}

\section{Introduction}

Let $U$ be a bounded open subset of $\mathbb{C}^n$. A weight $v$ on $U$ is a 
continuous, bounded, strictly positive real valued function on $U$. We will use 
${ \mathcal H}_v(U)$ to denote the Banach space of all holomorphic functions 
$f$ on $U$ 
which have the property that $\|f\|_v:=\sup_{z\in U}v(z)|f(z)|<\infty$ endowed 
with the norm  $\|\cdot\|_v$.  Consider all $f$ in ${\mathcal 
H}_v(U)$ with the property that $|f(z)|v(z)$ converges to $0$ as $z$
converges to the boundary of $U$ i.e. given $\epsilon>0$ there is a
compact subset, $K$, of $U$ such that $v(z)|f(z)|<\epsilon$ for $z$ in
$U\setminus K$. The set of all such functions is a closed subspace of 
${\mathcal H}_v(U)$ denoted by ${ \mathcal H}_{v_o}(U)$. We say that the weight 
$v$
on a balanced domain $U\subset \mathbb C^n$ is radial if $v(\lambda z)=v(z)$ 
for all $\lambda$ in $\mathbb{C}$ with $|\lambda|=1$ and all $z\in U$. 

In \cite{BR3,BR4,BR5} the authors characterised the surjective isometric 
isomorphisms of weighted
spaces of holomorphic functions. For radial weights on balanced open domains
in $\mathbb{C}$ they gave a complete characterisation. If $U$ and $V$ are 
bounded, balanced open subsets of $\mathbb{C}$ and $v\colon U\to\mathbb{R}$ and 
$w\colon
V\to\mathbb{R}$ are radial weights then every surjective isometric isomorphism 
$T\colon{ \mathcal H}_{v_o}(U)\to { \mathcal H}_{w_o}(V)$ has the form 
$$T(f)(z)=h_\phi(z) f\circ\phi(z)$$
for all $f$ in ${ \mathcal H}_{v_o}(U)$ and all $z$ in $V$
where $\phi\colon V\to U$ is a biholomorphic mapping and $h_\phi$ belongs to
${ \mathcal H}_{w_o}(V)$.  Since ${ \mathcal H}_{v_o}(U)$ is an M-ideal in 
${ \mathcal H}_{v}
(U)$, a theorem of Harmand and Lima, \cite{HL}, implies that every surjective 
isometric isomorphism $T\colon{ \mathcal H}_{v}(U)\to { \mathcal H}_{v}(U)$ also 
has the form 
$$T(f)(z)=h_\phi(z) f\circ\phi(z)$$
where, again, $\phi\colon U\to U$ is biholomorphic and $h_\phi$ belonging to
${ \mathcal H}_{v_o}(U)$.

In \cite{BLW} Bonet, Lindstr\"om and Wolf examined the isometric (not 
necessarily surjective) weighted composition operators between weighted
spaces of holomorphic functions giving both necessary and sufficient conditions
for a weighted composition operator to be an isometry. 
In this paper we shall examine the surjectivity of isometries between weighted 
spaces of holomorphic functions. We shall see that in many cases every 
isometry $T$ from ${ \mathcal H}_{v_o}(U)$ into ${ \mathcal H}_{v_o}(U)$ is 
automatically
surjective. Nevertheless, as we shall see, examples of non-surjective 
isometries from ${ \mathcal H}_{v_o}(U)$ into ${ \mathcal H}_{w_o}(U)$ and from 
${ \mathcal 
H}_{v}(U)$ into ${ \mathcal H}_{v}(U)$ do exist. Our examination of the weighted
spaces of holomorphic functions requires techniques 
from complex analysis, potential theory, topology and the geometry of Banach
spaces and prove that the surjectivity of
isometries is related to the separability and topological properties of a 
certain distinguished subspace of $V$, denoted by ${ \mathcal B}_v^T(U)$. 
For further reading on the isometric
theory of Banach spaces we refer the reader to \cite{AF,FJ1,FJ2}. For more
details on the geometric
theory and isometries of weighted spaces of holomorphic functions we refer the
reader to \cite{BiSu,BR1,BR2,BR3,BR4, BR5,BR6, Lu1} and 
\cite{Lu2}. 

We give  descriptions of the not necessarily surjective isometries 
complementing  former results for the surjective isometries
given in \cite{BR4}. Indeed, we show that in many cases this shows that all 
isometries are automatically surjective. Cima and Wogen \cite{CW} show
 the surjectivity of all 
isometries of the little Bloch space ${\mathcal B}_0$. As the space ${\mathcal 
B}_0$ is isometrically isomorphic to a particular weighted space of analytic 
functions, our results generalise those of Cima and Wogen.   
 In our generalisation of their results to Banach weighted spaces of 
holomorphic functions,
we are able to give a proof which includes both the Theorem and Proposition in 
\cite{CW}. Our proof is nontrivial and recovers their incomplete proofs (see 
the comments in Section 7). The key that let us afford the 
results with 
success is Theorem \ref{interior}, whose subtle proof uses Baire Category 
Theorem.

\section{Isometries of ${ \mathcal H}_{v_o}(U)$}\label{vo}

Let us begin this section with an example of non-surjective isometry
between weighted spaces of holomorphic functions. Let $\Delta$ 
denote the open unit disc in the complex plane.

\begin{example}\label{1to2}
Consider the weights $v(z)=1-|z|$ and $w(z)=1-|z|^2$ or the weights $v(z)=e^{
\frac{-1}{1-|z|}}$ and $w(z)=e^{\frac{-1}{1-|z|^2}}$ on the open unit disc $\Delta$.
Then, a routine calculation shows that, in both cases, $T(f)(z)=f(z^2)$ is a 
non-surjective isometry from ${ \mathcal H}_{v_o}(\Delta)$ into ${ \mathcal 
H}_{w_o}(\Delta)$.
\end{example}

In our classification of the surjective isometries of ${ \mathcal H}_{v_o}(U)$ a 
crucial role was fulfilled by a certain 
distinguished subspace of $U$, the $v$-boundary of $U$. 
In \cite{BR4} we showed that the set of extreme points of the closed
unit ball 
of ${ \mathcal H}_{v_o}(U)'$ is contained in $\{\lambda v(z)\delta_z: z\in U, 
|\lambda|=1\}$. The $v$-boundary of $U$ is defined as the set of all 
$z\in U$ such that $v(z)\delta_z$ is an extreme point of the unit ball of 
${ \mathcal H}_{v_o}(U)'$. Note that  $v(x)\delta_x$ is an extreme point of the 
unit 
ball of ${ \mathcal H}_{v_o}(U)'$ if and only if  $\lambda v(x)\delta_x$ is an 
extreme point for every $|\lambda|=1$. We denote the $v$-boundary of
$U$ by  ${ \mathcal B}_v(U)$. A weight $v$ is said to be complete if 
${\mathcal B}_v(U)=
U$. Each of the following weights on $B_{{\bf C}^n}$ is complete for $\alpha>0$,
$\beta\geqslant 1$
\begin{enumerate}
\item[(a)]  $v_{\alpha,\beta}(z)=(1-\|z\|^\beta)^\alpha$.
\item[(b)]  $\displaystyle w_{\alpha,\beta}(z)=e^{-\alpha\over 1-\|z\|^\beta}$.
\item[(c)]  $v(z)=(\log(2-\|z\|))^\alpha$.
\item[(d)]  $v(z)=(1-\log(1-\|z\|))^{-\alpha}$.
\end{enumerate}

When we consider (non necessarily surjective) isometries of 
${ \mathcal H}_{v_o}(U)$ into ${ \mathcal H}_{w_o}(V)$ we 
shall require a replacement for ${ \mathcal B}_w(V)$. As the extreme points of 
the 
closed unit ball of $T({ \mathcal H}_{v_o}(U))'$ are again contained in the set  
$\{\lambda w(z)\delta_z: z\in V, |\lambda|=1\}$ (see \cite[Lemma V.8.6 ]{DuSc})
we denote by ${ \mathcal B}^T
_w(V)$ the set of all $z\in V$ such that $w(z)\delta_z$ is an extreme point of 
the unit ball of $T({ \mathcal H}_{v_o}(U))'$. Note that different $z$ in 
${ \mathcal B}^T_w(V)$ may lead to the same extreme point in  the closed unit 
ball of $T({ \mathcal H}_{v_o}(U))'$. 

The following result is a Banach-Stone type theorem for isometries between 
weighted 
spaces of holomorphic functions. Its proof is heavily modelled on that given by 
Cima 
and Wogen, \cite{CW}, which characterises  the isometries 
of the little Bloch space with the weight 
$1-|z|^2$ replaced with an arbitrary radial weight. We include its proof for 
completeness and future reference. See also \cite[Theorem~3.1]
{AF}.

\begin{theorem}\label{stone}
Let  $V$ be an open subset of $\mathbb{C}$. Let $v\colon \Delta\to \mathbb 
{R}^+$, $w\colon V\to \mathbb{R}^+$ be weights with $v$ radial and converging 
to $0$ on the boundary of $\Delta$. Let $T\colon { \mathcal 
H}_{v_o}(\Delta)\to { \mathcal H}_{w_o}(V)$ be 
an isometry. Then there is a holomorphic  function $\phi\colon V
\to \Delta$ and $h_\phi$ in ${ \mathcal H}_{w_o}(V)$ such that
$$T(f)(z)=h_\phi(z) f\circ\phi(z)$$
for all $f$ in ${ \mathcal H}_{v_o}(\Delta)$ and all $z$ in $V$.
\end{theorem}

{\sc Proof:} We assume without loss of generality that $v(z)\leqslant 1$ for
all $z$ in $\Delta$. Consider the surjective isometry $T\colon { \mathcal 
H}_{v_o}(\Delta)\to T({ \mathcal H}_{v_o}(\Delta))$. Then $T'$, the transpose of 
$T$, maps $T({ \mathcal H}_{v_o}(\Delta))'$ 
isometrically onto $({ \mathcal H}_{v_o}(\Delta))'$. Hence it maps the extreme 
points
of the unit ball of $T({ \mathcal H}_{v_o}(\Delta))'$ bijectively onto the set 
of extreme 
points of the unit ball of $({ \mathcal H}_{v_o}(\Delta))'$. This induces a 
surjective function  
$\phi_1\colon { \mathcal B}_{w}^T(V)\to { \mathcal B}_v(\Delta)$ and a function 
 $\alpha\colon  { \mathcal B}_w^T(V)\to \mathbb C$ with $|\alpha(\cdot)|=1$
 so that $T'(w(z)\delta_z)=\alpha(z)v(\phi_1(z))\delta_{\phi_1(z)}$ 
for all $z$ in ${ \mathcal B}^T_w(V)$. Let $f_o\equiv1$. Then 
$$w(z){T}(f_o)(z)=\alpha(z)v\circ \phi_1(z)$$
for all $z$ in ${ \mathcal B}_w^T(V)$ 
or 
\begin{equation}\label{mod}
T(f_o)(z)=\alpha(z)\frac{v\circ\phi_1(z)}{w(z)}
\end{equation}
for all $z$ in ${ \mathcal B}_w^T(V)$. Note that $T(f_o)(z)\not=0$ for $z$ in 
${ \mathcal B}_w^T(V)$. Similarly taking $f_1(z)=z$ for $z$ in $V$ we get 
$$T(f_1)(z)=\alpha(z)\frac{v\circ\phi_1(z)}{w(z)}\phi_1(z)$$
for all $z$ in ${ \mathcal B}_w^T(V)$. This gives us that 
$$\phi_1(z)=\frac{T(f_1)(z)}{T(f_o)(z)}$$
for $z$ in ${ \mathcal B}_w^T(V)$. We note that the right-hand side of the 
above 
equation is defined and holomorphic for all $z$ in $V\setminus T(f_o)^{-1}(0)$. 
As $v$ is radial ${ \mathcal B}_w^T(V)$ is uncountable (see 
\cite[Lemma 5]{BR1}) and therefore has an  accumulation point in $V$. In 
particular, this means that we can extend $\phi_1$ to a holomorphic
function $\phi_2$ on $V\setminus T(f_o)^{-1}(0)$ by setting it equal to
$\displaystyle{\frac{T(f_1)(z)}{T(f_o)(z)}}$. Since $\Delta$ 
is bounded and $v$ converges to $0$ on the boundary of $\Delta$ it follows that
$\mathcal{H}_{v_o}(\Delta)$ contains all polynomials. Next we consider the 
function $f_k(z)=z^k$. We get that 
$$T(f_k)(z)=T(f_o)(z)\phi_1(z)^k$$
for all $z$ in ${ \mathcal B}_w^T(V)$. The Identity Principle gives that 
$$T(f_k)(z)=T(f_o)(z)\phi_2(z)^k$$
and thus
$$
w(z)|T(f_k)(z)|=w(z)|T(f_o)(z)||\phi_2(z)|^k
$$
for all $z$ in $V\setminus T(f_o)^{-1}(0)$. Taking $k$th roots and letting $k$
tends to infinity we observe that $\phi_2$ is bounded on 
$V\setminus T(f_o)^{-1}(0)$. This means that we can extend $\phi_2$ analytically
to a 
holomorphic function on $V$ which we denote by $\phi$. We claim as $\Delta$ 
is convex  then $\phi(V)\subseteq \Delta$. First we show that $\phi(V)\subseteq 
\overline{\Delta}$. Suppose this is not the case. Then we can 
choose a continuous linear functional, $l$, on $\mathbb{C}$ with $\|l\|_\Delta
\leqslant 1$
 so that $|l(\phi(z_0))|>1$ for some $z_0\in V$. Continuity 
allows us to suppose in addition that $T(f_o)(z_0)\not=0$. We 
have that 
$$
T(l^k)(z_0)=T(f_o)(z_0) l(\phi(z_0))^k,
$$
and thus
$$
|T(l^k)(z_0)|w(z_0)=|T(f_o)(z_0)|\ w(z_0)|l( \phi(z_0))|^k,
$$
for all $k\in {\bf N}$. Letting $k$ tend to infinity gives a contradiction and
therefore $|l(\phi(z))|\leqslant 1$ for all linear $l$ with
$\|l\|_\Delta\leqslant 1$. Therefore $\phi(V)\subseteq \overline{\Delta}$. The 
Open Mapping Theorem implies that  $\phi(V)\subseteq \Delta$. A final 
application of the Identity Principle gives that 
$$T(f)(z)=h_\phi(z) f\circ\phi(z),$$
where $h_\phi:=T(f_o)$, for all $f$ in ${ \mathcal H}_{v_o}(\Delta)$ and all $z$ 
in $V$.

We note that the $\phi$ of the above theorem maps ${ \mathcal B}_w^T(V)$ onto 
${ \mathcal B}_v(\Delta)$.

We can replace the condition that $v$ is radial with the 
condition that $v$ is complete and obtain the following theorem.

\begin{theorem}\label{stonec}
Let  $V$ be an open subset of $\mathbb{C}$  and $v\colon \Delta\to 
\mathbb {R}^+$, $w\colon V\to \mathbb{R}^+$ be weights with $v$ complete. Let 
$T\colon { \mathcal H}_{v_o}(\Delta)\to { \mathcal H}_{w_o}(V)$ be 
an isometry. Then there is an analytic surjection $\phi\colon V\to \Delta$ and 
$h_\phi$ in ${ \mathcal H}_{w_o}(V)$ such that
$$T(f)(z)=h_\phi(z) f\circ\phi(z)$$
for all $f$ in ${ \mathcal H}_{v_o}(\Delta)$ and all $z$ in $V$.
\end{theorem}

We note that from (\ref{mod}) and by continuity we have that
$$|h_\phi(z)|=\frac{v\circ\phi(z)}{w(z)}$$
for $z$ in $\overline{{\mathcal B}_w^T(V)}\cap V$.

\section{Separation condition and Surjectivity of Isometries}

Let $X$ be a Hausdorff locally compact space. Denote by $C_o(X)$ the space of 
continuous $\mathbb{C}$-valued functions on $X$ which vanish at infinity. 
According 
to Araujo and Font, \cite{AF}, a subspace $A$ of $C_o(X)$ is strongly separating
if for each $x,y$ in $X$ with $x\not=y$ there is $f$ in $A$ with $|f(x)|\not=
|f(y)|$. 
In \cite{AF} Araujo and Font showed that this separation condition on the 
range
of an isometry allows  significantly stronger Banach--Stone type theorems.
Our main result of this section shows that under certain conditions strong 
separation is a necessary and sufficient condition for the surjectivity of
isometries between spaces of type ${ \mathcal H}_{v_o}(U)$. Let $V$ be an open 
subset
of $\mathbb{C}^n$. We use $\overline{V}^*$ to
denote the one-point compactification of $V$. We shall say that a subspace $A$
of ${\mathcal H}_w(V)$ {\it strongly $w$-separates the points of 
$V$\/}
if for each $x,y$ in $V$ with $x\not=y$ there is $f$ in $A$ with $w(x)|f(x)|
\not=w(y)|f(y)|$.

\begin{theorem}\label{separ}
Let $V$ be a connected open subset of $\mathbb{C}$. Let $v\colon \Delta\to 
\mathbb{R}^+$, $w\colon V\to \mathbb{R}^+$ be continuous weights with $v$ 
complete and converging to $0$ on the boundary of $\Delta$ and $w$ be such that 
${ \mathcal H}_{w_o}(V)$ contains all polynomials of degree $1$. Let $T\colon 
{ \mathcal H}_{v_o}(\Delta)\to { \mathcal H}_{w_o}(V)$ be an isometry. Then the 
following are equivalent
\begin{enumerate}
\item[(a)] $T$ is surjective,
\item[(b)] $T({ \mathcal H}_{v_o}(\Delta))$ strongly $w$-separates the points of 
$V$,
\item[(c)] $T({ \mathcal H}_{v_o}(\Delta))$ contains all polynomials of degree 
$1$.
\end{enumerate}
\end{theorem}

{\sc Proof:} If $T({ \mathcal H}_{v_o}(\Delta))$ contains all 
polynomials of 
degree $1$ then $T({ \mathcal H}_{v_o}(\Delta))$ will  strongly $w$-separates 
the points of $\overline{V}^*$. 

By Theorem \ref{stonec} there is an analytic surjection $\phi\colon V\to\Delta$ 
and 
$h_\phi$ in ${ \mathcal H}_{w_o}(V)$ such that
$$T(f)(z)=h_\phi(z) f\circ\phi(z)$$
for all $f$ in ${ \mathcal H}_{v_o}(\Delta)$ and all $z$ in $V$. Moreover for 
all $z$ in ${ \mathcal B}_w^T(V)$ we have that 
$$|h_\phi(z)|=\frac{v\circ\phi(z)}{w(z)}.$$

Now suppose that $T({ \mathcal H}_{v_o}(\Delta))$ strongly $w$-separates the 
points of 
$V$. Consider the isometry $T^{-1}\colon T({ \mathcal H}_{v_o}(\Delta))\to 
{ \mathcal H}_{v_o}(\Delta)$. Then $(T^{-1})'$ is an isometry of ${ \mathcal 
H}_{v_o}(\Delta)'$ onto $(T({ \mathcal 
H}_{v_o}(\Delta))'$. As such $(T^{-1})'$ maps the extreme points of the unit ball
of 
${ \mathcal H}_{v_o}(\Delta)'$ onto the extreme points of the unit ball of 
$(T({ \mathcal H}_{v_o}(\Delta))'$. Hence for each $x$ in $\Delta$ there is 
$\mu$ in 
$\mathbb{C}$ with $|\mu|=1$ and $y$ in ${ \mathcal B}_w^T(V)$ so that
$$(T^{-1})'(v(x)\delta_x)=\mu w(y)\delta_y.$$ 
Since $T({ \mathcal H}_{v_o}(\Delta))$ strongly $w$-separates the points of $
V$ each extreme point of the unit ball of $(T({ \mathcal H}_{v_o}(\Delta))'$ 
determines a unique $y$ in ${ \mathcal B}_w^T(V)$ and hence there is a function 
$\beta\colon U\to \mathbb C$ with $|\beta(\cdot)|=1$, and a function 
$\psi\colon U\to { \mathcal B}_{w}^T(V)$ such that $(T^{-1})'(v(x)\delta_x)=
\beta(x) 
w(\psi(x))\delta_{\psi(x)}$ for all $x$ in $\Delta$. (See also 
\cite[Theorem~3.1]{AF}.)
It now follows that there is a function $h_{\psi}\colon \Delta\to\mathbb{C}$ 
such that 
$$T^{-1}(f)(z)=h_{\psi}(z)f\circ\psi(z)$$
for all $f\in T({ \mathcal H}_{v_o}(\Delta))$ and all $z$ in $\Delta$. 
Moreover, we also have that
$$|h_{\psi}(z)|=\frac{w\circ \psi(z)}{v(z)}$$
for all $z$ in $\Delta$. 

We claim that $\psi$ is continuous. To see this consider a sequence $(z_n)_n$ 
in $\Delta$ converging to some point $z_0$ in $\Delta$. Since $|\beta(z_{n})|=1$ 
and 
$\psi(z_{n})$ is in the compact subset $\overline {{ \mathcal B}^T_{w}(V)}$ of
$\overline{V}^*$, for all
$k$, we can assume without loss of generality that there is a subsequence
$(z_{n_k})_k$ of $(z_n)_n$ so that $(\beta(z_{n_k}))_k$ converges to some
$\beta_0$ and $(\psi(z_{n_k}))_k$ converges to some $u_0\in \overline 
{{ \mathcal B}^T_{w}(V)}$. Then,
$(T^{-1})'(v(z_{n_k})\delta_{z_{n_k}})=\beta(z_{n_k})w(\psi(z_{n_k}))\delta_{
\psi(z_{n_k})}$ converges weak$^*$ to
$(T^{-1})'(v(z_0)\delta_{z_0})=\beta(z_0)w(\psi(z_0))\delta_{\psi(z_0)}$ and to
$\beta_0w(u_0)\delta_{u_0}$. Since $T({ \mathcal H}_{v_o}(\Delta))$ separates 
points of 
$V$ we have that $u_0\in V$ and  
$u_0=\psi(z_0)$. Hence $\psi(z_{n_k})$ converges to $\psi(z_0)$. Applying the
above argument to any subsequence of $(z_n)_n$ we get that $\psi$ is 
continuous.

We next observe that $\phi\circ\psi(z)=z$ for $z$ in $\Delta$. 
To see this we note that for every $f$ in ${ \mathcal H}_{v_o}(\Delta)$ we have 
$$f(z)=T^{-1}(T(f))(z)=h_\psi(z)h_\phi(\psi(z))f\circ\phi\circ\psi(z)$$
for all $z$ in $\Delta$. Taking $f\equiv 1$ we get $h_\psi(z)h_\phi(\psi(z))=1$ 
for 
all $z$ in $\Delta$ and this gives that $\phi\circ\psi={\rm Id}_\Delta$. 

We also note that 
$$g(z)=T(T^{-1}(g))(z)=h_\phi(z)h_\psi(\phi(z))g\circ\psi\circ\phi(z)$$
for all $g\in T({ \mathcal H}_{v_o}(\Delta))$, $z$ in $V$ which gives that
$$w(z)|g(z)|=w(\psi\circ\phi(z))|g\circ\psi\circ\phi(z)|$$
for all $g$ in $T({ \mathcal H}_{v_o}(\Delta))$, $z$ in ${ \mathcal B}^T_w(V)$. 
Since
$T({ \mathcal H}_{v_o}(\Delta))$ strongly $w$-separates the points of $V$ we get 
that
$\psi\circ\phi(z)=z$ for all $z$ in ${ \mathcal B}^T_w(V)$. 

As $\phi\circ\psi=Id_\Delta$ we have that $\psi$ is injective. The Invariance of 
Domains implies that ${ \mathcal B}_w^T(V)=\psi(\Delta)$ is open in 
$\mathbb{C}$ and 
hence $\psi=\left(\phi|_{{ \mathcal B}_w(V)}\right)^{-1}$ and $h_\psi$ are therefore  
holomorphic on $\Delta$. 

Let $g$ belong to ${ \mathcal H}_{w_o}(V)$. We define $f\colon U\to \mathbb{C}$ 
by 
$$f(z)=\frac{g\circ\psi(z)}{h_\phi(\psi(z))}.$$ 
Note that since $\psi(\Delta)={ \mathcal B}_{w}^T(V)$ we have that $|h_\phi
(\psi(z))|=\frac{v(z)}{w(\psi(z))}\not=0$ for all $z$ in $\Delta$. This, in 
particular, will mean that $f$ is well defined and holomorphic on $\Delta$. 

We claim that $f\in { \mathcal H}_{v_o}(\Delta)$. To see this let $\epsilon>0$ 
and 
take a compact set $K\subset V$ such that $w(x)|g(x)|<\epsilon$ for all 
$x\in V\setminus K$. By continuity $\phi(K)$ is a compact set in $\Delta$. Since
$\phi\circ\psi(z)=z$ for all $z$ in $\Delta$, for every $z\in \Delta\setminus 
\phi(K)$ it follows that $\psi(z) \in V\setminus K$. Then
$$
v(z)|f(z)|=v(z)\frac{|g\circ \psi(z)|}{|h_\phi(\psi(z))|}=w(\psi(z))|g\circ 
\psi(z)|<\epsilon
$$
for $z\in \Delta\setminus \phi(K)$.

Finally for $z=\psi\circ\phi(z)$ in ${ \mathcal B}_{w}^T(V)$ we have that
$$
T(f)(z)=T\left(\frac{g\circ\psi}{h_\phi\circ\psi}\right)(z)
=h_\phi(z)\frac{g\circ\psi(\phi(z))}{h_\phi(\psi(\phi(z))}
=g(z).
$$
As $V$ is connected and ${ \mathcal B}_{w}^T(V)$ is open the Identity Principle 
implies that $T(f)(z)=g(z)$ for  all $z$ in $V$ and hence $T$ is surjective.

Clearly if $T$ is surjective then $T({ \mathcal H}_{v_o}(\Delta))$ will contain 
all the linear functional and the proof is complete.

We observe that the condition that $V$ is connected is necessary. To see this
let $U_1$ and $U_2$ be two disjoint non-empty open subset of $\mathbb{C}$ and 
$v_1\colon U_1\to\mathbb{R}^+$ and $v_2\colon U_2\to \mathbb{R}^+$ be complete 
continuous weights such that ${ \mathcal H}_{{(v_2)}_o}(U)\not=\{0\}$. Let 
$V=U_1\cup U_2$ and $w\colon V\to \mathbb{R}$ be given by 
$$w(z)=\begin{cases} v_1(z)&\text {if $z\in U_1$,}\\
                              v_2(z) & \text{if $z\in U_2$.}\\
\end{cases}
$$
Then the mapping $T\colon { \mathcal H}_{{v_1}_o}(U_1)\to { \mathcal H}_{w_o}(V)$ 
given by
$$T(f)(z)=\begin{cases} f(z) &\text {if $z\in U_1$,}\\
                     0 & \text{if $z\in U_2$.}\\
\end{cases}
$$
is an non-surjective isometry of ${ \mathcal H}_{v_o}(U_1)$ into 
${ \mathcal H}_{w_o}(V)$. 

We also note that the condition of separability of Theorem~\ref{separ} is not 
satisfied by the isometry in Example~\ref{1to2} as $T(f)(z)=f(z^2)$ can never 
separate $z$ and $-z$. 
To see that this isometry cannot satisfy condition (c) of Theorem~\ref{separ}
we observe that for any $f$ in ${ \mathcal H}_{v_o}(\Delta)$ the restriction of 
$T(f)$ to $\mathbb{R}$ is an even function. Hence it cannot be equal to any 
polynomial of degree one. 

\section{Topological Structure of ${ \mathcal B}_w^T(V)$}

In subsequent sections we will wish to consider isometries between specific 
spaces of holomorphic functions. Previous sections have shown that vital
information of the structure of isometries is contained in the subspace 
${\mathcal 
B}_w^T(V)$ of $V$. The results in this section will cast light on the 
topological nature of ${ \mathcal B}_v^T(V)$. 

\begin{theorem}\label{interior}
 Let $U$ and $V$ be bounded open subsets of $\mathbb{C}$ and $\phi\colon V\to
U$ be a surjective analytic function. Let $B$ be a subset of $V$ such that
$\phi(B)=U$.  Then for every $a$ in $U$ and every $r>0$, $\overline{
B\cap \phi^{-1}(B(a,r))}$ has non-empty interior. In particular, $\overline B$ 
has non-empty interior.
\end{theorem}

{\sc Proof:} Let us suppose that we can find $a$ in $U$ and $r>0$ so
that the interior of $\overline{B\cap \phi^{-1}(B(a,r))}$ is empty. As 
$\phi\colon V\to 
U$ is analytic and non-constant we have that $\phi$ is locally injective on 
all of $V$ with the possible exception of a sequence of points, $(z_n)_n$ in 
$V$, which converges to the boundary of $V$. For each $k$ in $\mathbb{N}$ let 
$$
E_k=\left(\overline{B\cap \phi^{-1}(B(a,r))}\cap \left\{z\in V:{\rm dist}(z,
\partial V)\geqslant \frac{1}{k}\right\}\right)\setminus\bigcup_{n}B\left(z_n,
\frac{1}{2^{k+2}}\right).
$$
Then $\{E_k\}_k$ is a family of compact subsets of $\overline{B\cap \phi^{-1}(B
(a,r))}\cap V$ such 
that $\bigcup_{k=1}^\infty E_k=\left(\overline{B\cap \phi^{-1}(B(a,r))} 
\setminus\{(z_n)_n\}\right)\cap V$. 

For each $z$ in $\left(\overline{B\cap \phi^{-1}(B(a,r))}\setminus\{(z_n)_n\}
\right)\cap V$ we choose $\delta_z>0$ so that $\phi$ is injective on $B(z,
\delta_z)$ and let $W_z=\overline{B(z,\delta_z/2)}$. As each $E_k$ is compact 
for each $k$ in $\mathbb{N}$ we
can choose $z_1^k,z_2^k,\ldots,z_p^k$ in $\left(\overline{B\cap \phi^{-1}
(B(a,r))}\setminus\{(z_n)_n\}\right)\cap V$ so that $E_k\subset 
\bigcup_{j=1}^pW_{z_j^k}$. Since $\overline{B\cap \phi^{-1}(B(a,r))}$ has empty 
interior, $E_k\cap W_{z_j^k}$, $j=1,\ldots, p$, will have empty interior. As the 
restriction of $\phi$ to $W_{z_j^k}$ is a homeomorphism onto its image it 
follows that $\phi(E_k\cap W_{z_j^k})$ is a compact and hence closed set which 
will also have empty interior. Therefore, by the 
Baire Category Theorem, $\phi(E_k)$ has empty interior. It is easy to see that
$$
U\cap B(a,r)\subseteq \phi \left(\overline{B\cap\phi^{-1}(B(a,r))}\cap V 
\right) \subseteq U\cap\overline{B(a,r)}.
$$
Hence 
$$
\bigcup_{k=1}^\infty \phi(E_k)\cup \overline{\{\phi(z_n):n\in {\mathbb N}\}}\cup
\partial B(a,r)= 
\overline{B(a,r)}\cup \overline{\{\phi(z_n):n\in {\mathbb N}\}}.
$$
This contradicts the Baire Category Theorem and hence we have that 
$\overline{B\cap \phi^{-1}(B(a,r))}$ has non-empty interior.

We observe that Theorem~\ref{interior} fails if we replace $\overline{B\cap 
\phi^{-1}(B(a,r))}$ with $B\cap \phi^{-1}(B(a,r))$. To see this consider the
surjective analytic function $\phi\colon\Delta\to\Delta$ given by $\phi(z)=
z^2$. We consider the subset $B$ of $\Delta$ given by
$$
B=\{re^{i\theta}:r\in \mathbb{Q},0\leqslant\theta\leqslant \pi\}\cup \{re^{i\theta}:r\in 
\mathbb{R}\setminus\mathbb{Q},\pi\leqslant\theta\leqslant 2\pi\}.
$$
Then $\phi(B)=\Delta$ yet $B$ does not have an interior point. 

Taking $B={ \mathcal B}_w^T(V)$ we get the following result.

\begin{corollary}\label{int2}
 Let  $V$ be a bounded open subset of $\mathbb{C}$ and  Let 
$v\colon \Delta\to \mathbb {R}^+$, $w\colon V\to \mathbb{R}^+$ be 
weights with $v$ complete. Let $T\colon { \mathcal H}_{v_o}(\Delta)\to 
{ \mathcal H}_{w_o}(V)$ be 
an isometry and ${ \mathcal B}_w^T(V)$, $\phi\colon V\to \Delta$ be as in 
Theorem~\ref{stonec}. Then for every $a$ in $\Delta$ and every $r>0$, 
$\overline{
{ \mathcal B}_w^T(V)\cap \phi^{-1}(B(a,r))}$ has non-empty interior.
\end{corollary}

Let us see that we cannot replace the assumption that $v$ is complete in 
Corollary~\ref{int2}
with the assumption that $v$ is radial. Consider the weight $v\colon \Delta \to 
\mathbb R$ given by $v(z)=1-c(|z|)$, $z\in \Delta$, where $c\colon [0,1]\to
[0,1]$ is the Cantor function (see \cite[Problem 1.5.20]{MP}). Then 
${ \mathcal B}_v(\Delta)$
consists of a countable collection of circles centred at $0$ and as such has 
empty interior.
If we consider the identity mapping $T\colon { \mathcal H}_{v_o}(\Delta)\to 
{ \mathcal H}_{v_o}(\Delta)$, $T(f)=f$, then $\overline{{ \mathcal B}^T_v
(\Delta)}=
\overline{{ \mathcal B}_v(\Delta)}={ \mathcal B}_v
(\Delta)\cup\{0\}\cup\{z:|z|=1\}$ and therefore also has empty interior.

\section{Automatic Surjectivity of Isometries}

In this section we will show that for the weights $v(z)=1-|z|^\beta$, 
$\displaystyle{v(z)=e^{\frac{-1}{1-|z|^\beta}}}$ $\beta\geqslant 1$ and the weight
$v(z)=(1-\log(1-|z|))^\beta$, $\beta<0$, on the unit disc $\Delta$ we have 
automatic surjectivity of isometries from ${ \mathcal H}_{v_o}(\Delta)$ into 
${ \mathcal H}_{v_o}(\Delta)$. In all three cases the general approach is the 
same.
Using Theorem~\ref{stone} we know that each isometry $T$ has the form $T(f)(z)
=h_\phi(z)f\circ\phi(z)$. With the exception of the weight $v(z)=1-|z|^2$ we then
show that $\phi(0)=0$. From this we can then prove that $\phi$ is an 
automorphism of the disc and hence that $T$ is surjective. While our general
strategy is the same in all three cases we are forced not only to use different
arguments for each individual weight but also different arguments for different
values of $\beta$ in each of the first two cases.   In 
\cite[Theorems~13, 15 and 16]{BR4} the  
surjective isometries of ${ \mathcal H}_{v_o}(\Delta)$ are completely described.

\begin{lemma}
Let  $v\colon \Delta\to \mathbb {R}^+$ be a radial or 
complete weight and $T\colon { \mathcal H}_{v_o}(\Delta)\to { \mathcal H}_{v_o}
(\Delta)$ be an isometry. If there exists an automorphism $\phi:\Delta\to 
\Delta$ such that $T(f)(z)=\phi'(z) f\circ \phi(z)$ for all $z\in \Delta$ 
then, $T$ is surjective.
\end{lemma}

{\sc Proof:}
Define $S\colon { \mathcal H}_{v_o}(\Delta)\to {\mathcal 
H}(\Delta)$ by $S(g)(z):=(\phi^{-1})'(z)g\circ\phi^{-1}(z)$. It is easily checked 
that $S$ maps ${ \mathcal H}_{v_o}(\Delta)$ into ${ \mathcal H}_{v_o}(\Delta)$
and that $S=T^{-1}$ proving that $T$ is surjective.

 Let us start with the weight 
$v(z)=1-|z|^\beta$. 

\begin{theorem}\label{beta1}
Let $\beta\geqslant 1$ and $v\colon\Delta\to\Delta$ be given by 
$v(z)=1-|z|^\beta$. 
Let $T\colon { \mathcal H}_{v_o}(\Delta)\to { \mathcal H}_{v_o}(\Delta)$ be an 
isometry.
\begin{enumerate}
\item[(a)] If $\beta=2$ then there exists an automorphism $\varphi:\Delta\to 
\Delta$
such that $T(f)(z)=\varphi'(z) f\circ \varphi(z)$ for all $z\in \Delta$.
\item[(b)] If $\beta\neq 2$ then there exist $ \theta\in \mathbb R$ and a 
complex number $\alpha$, $|\alpha|=1$, such that $T(f)(z)=\alpha 
f(ze^{i\theta})$ for all $z\in \Delta$.
\end{enumerate}
In particular $T$  
is surjective.
\end{theorem}

{\sc Proof:} By Theorem~\ref{stone} we know that there is an analytic
surjection $\phi\colon\Delta\to\Delta$ and $h_\phi$ in ${ \mathcal H}_{v_o}
(\Delta)$ such that 
$$T(f)(z)=h_\phi(z) f\circ\phi(z)$$
for all $f$ in ${ \mathcal H}_{v_o}(\Delta)$ and all $z$ in $\Delta$. 
Theorem~\ref{interior} tells us that the interior of $\overline{{ \mathcal 
B}_v^T(\Delta)}$ is non-empty. Moreover, for points of this set we have that
$$|h_\phi(z)|=\frac{v\circ\phi(z)}{v(z)}.$$ As $h_\phi$ is analytic and non-zero
on the interior of $\overline{{ \mathcal B}_v^T(\Delta)}$ we have that $\log 
|h_\phi(z)|$ is harmonic on $\overline{{ \mathcal B}_v^T(\Delta)}$. Hence we 
have that
$$\Delta\log(1-|\phi(z)|^\beta)=\Delta\log(1-|z|^{\beta})$$ 
or that
\begin{equation}\label{harlog}
\frac{|\phi(z)|^{\beta-2}|\phi'(z)|^2}{(1-|\phi(z)|^\beta)^2}=\frac{|z|^{\beta-2}}{
(1-|z|^{\beta})^2}
\end{equation}
for all $z$ in the interior of $\overline{{ \mathcal B}_v^T(\Delta)}$.

We consider four cases depending on the value of $\beta$. In 
each of these cases we will show that $\phi$ is an
automorphism of the disc which implies that $T$ is surjective. 

When $\beta=2$, Equation~\eqref{harlog} becomes
$${|\phi'(z)|^2}=\frac{(1-|\phi(z)|^2)^2}{(1-|z|^2)^2}$$
for $z$ in the interior of $\overline{{ \mathcal B}_v^T(\Delta)}$.
Applying the Schwarz-Pick Lemma we see that $\phi$ must be an automorphism of 
the disc. 

We now consider the other cases. For each $n\in\mathbb{N}$ the
set $\overline{{ \mathcal B}_v^T(\Delta)\cap \phi^{-1}(B(0,\frac{1}{n}))}$ has 
non-empty interior. Hence for each $n$ in $\mathbb{N}$ we can choose $a_n$ in 
the interior of $\overline{{ \mathcal B}_v^T(\Delta)\cap\phi^{-1}
(B(0,\frac{1}{n}))}$, $a_n\not=0$. Then $(\phi(a_n))_n$ is a null 
sequence in $\Delta$. Since $\overline{\Delta}$
is compact $(a_n)_n$ has a subsequence, which we also denote by $(a_n)_n$, that
converges to some point $a$ of $\overline{\Delta}$. We claim that $a$ is 
actually in $\Delta$. To see this suppose that $a$ belongs to $\partial\Delta$. 
Since each $a_n$ belongs to $\overline{{ \mathcal B}_v^T(\Delta)}$ and $h_\phi$ 
is in ${ \mathcal H}_{v_o}(\Delta)$ we have that
$$
\lim_{n\to\infty}(1-|a_n|^\beta) |h_\phi(a_n)|=\lim_{n\to\infty}(1-|\phi(a_n)|^\beta)=0
$$
contradicting the fact that $(\phi(a_n))_n$ is a null sequence. By continuity 
of $\phi$ we have that $\phi(a)=0$. Our aim is to prove that $a=0$.
 
Let us  consider the case $\beta>2$.
As $(a_n)_n$ converges to $a$ in $\Delta$ and 
$$\frac{|\phi(a_n)|^{\beta-2}|\phi'(a_n)|^2}{(1-|\phi(a_n)|^\beta)^2}=\frac{|a_n|^
{\beta-2}}{(1-|a_n|^{\beta})^2}$$
for each $n$ in $\mathbb{N}$ it follows that $(a_n)_n$ must be a null sequence.
By continuity of $\phi$ it follows that $\phi(0)=0$. We can 
apply the Schwarz Lemma to get that $\displaystyle{\frac{(1-|\phi(z)|^\beta)^2}
{(1-|z|^{\beta})^2}}\geqslant 1$ for all $z$ in $\Delta$. 

On the other hand, rewriting Equation \eqref{harlog} we have that
\begin{equation*}\label{powers}
{|\phi'(z)|^2}\left(\frac{|\phi(z)|}{|z|}\right)^{\beta-2} =\frac{(1-|\phi(z)|
^\beta)^2}{(1-|z|^{\beta})^2}
\end{equation*}
for $z\not=0$ 
in ${ \mathcal B}^T_v(\Delta)$.
Letting $z$ tend to $0$ we get that $|\phi'(0)|^\beta\geqslant 1$ and 
applying the Schwarz Lemma again we see that $|\phi'(0)|=1$ and that therefore 
$\phi$ is an automorphism of the disc. 

We now consider the case where $1<\beta<2$. Let $\gamma=\frac{2}{2-\beta}$. 
We rewrite equation~\eqref{harlog} as 
$$\left(\frac{1-|\phi(z)|^\beta}{1-|z|^\beta}\right)^\gamma|\phi(z)|=
|z||\phi'(z)|^\gamma$$
or as 
$$|h_\phi(z)|^\gamma|\phi(z)|=|z||\phi'(z)|^\gamma$$
for all $z$ in the interior of $\overline{{ \mathcal B}_v^T(\Delta)}$. 
Let us see that this equality can be extended to the whole of $\Delta$. 
Taking logs of both sides, we get that
$$
\gamma \log|h_\phi(z)|+\log|\phi(z)|=\log|z|+\gamma \log|\phi'(z)|
$$
for $z\in { \mathcal B}_v^T(\Delta)$. As $\gamma \log|h_\phi(z)|+\log|\phi(z)|$ 
and $\log|z|+\gamma \log|\phi'(z)|$ are harmonic on $\Delta\setminus 
( h_\phi^{-1}(0) \cup \phi^{-1}(0) \cup (\phi')^{-1}(0) \cup \{0\})$,
the Identity Principle implies that 
$$
\gamma \log|h_\phi(z)|+\log|\phi(z)|=\log|z|+\gamma \log|\phi'(z)|
$$
for $z\in\Delta\setminus 
( h_\phi^{-1}(0) \cup \phi^{-1}(0) \cup (\phi')^{-1}(0) \cup \{0\})$ 
and by continuity hence 
$$
|h_\phi(z)|^\gamma|\phi(z)|=|z||\phi'(z)|^\gamma
$$
 for all $z$  in $\Delta$. We claim that $\phi(0)=0$. Suposse this is not the 
case.
Then $h_\phi$ has a zero of order $k$ at $0$. If we let $l$ be the order of
 $\phi'$ at $0$ then we get $(k-l)\gamma=1$, which is impossible as 
$ \gamma>2$ and $k$ and $l$ are non negative integers.


Suppose that 
$a\not=0$ and that $\phi$ has a zero of degree $m$ at $a$. Let us 
first consider the case when $m=1$. Then, as $|h_\phi(z)|^\gamma|\phi(z)|=| z||
\phi'(z)|^\gamma$ for all 
$z$ in $\Delta$ we get that the left--hand--side has a zero of order $1$ at 
$a$ while the right--hand--side is non-zero at $a$. Now suppose that $m>1$. 
Since $|h_\phi(z)|^\gamma|\phi(z)|=| z||\phi'(z)|^\gamma$ for all $z$ in $\Delta$ 
we see
that $\displaystyle{m=\frac{\gamma}{\gamma-1}}$. But this now implies that 
$\displaystyle{\gamma=\frac{m}{m-1}}$ and again $\gamma=1$.  
This is impossible as we have assumed that $1<\beta<2$ and so,
as $0$ is the only zero of $\phi$ in $U$, $a=0$.

For each $n$ in $\mathbb{N}$ we have that 
$$
{|\phi'(a_n)|^2}\left|\frac{a_n}{\phi(a_n)}\right|^{2-\beta} =\frac{(1-
|\phi(a_n)|^\beta)^2}{(1-|a_n|^\beta)^2}.
$$
Letting $n$ tend to $\infty$ we get that $|\phi'(0)|\geqslant 1$. Applying the 
Schwarz Lemma again we see $\phi$ is an automorphism of the disc. 

Finally we consider the case when $\beta=1$. In this case we see that $\phi$
satisfies the  equation 
$$
\left|\frac{(1-|\phi(z)|)^2}{(1-|z|)^2}\right||\phi(z)|=|z||\phi'(z)^2|
$$
or as 
$$|h_\phi(z)^2\phi(z)|=|z\phi'(z)^2|$$
for all $z$ in the interior of $\overline{{ \mathcal B}_v^T(\Delta)}$. As in 
the case where $1<\beta<2$ we see that 
$$
|h_\phi(z)^2\phi(z)|=|z\phi'(z)^2|
$$
for all $z$ in $\Delta$. The Open Mapping Theorem allows us to find 
 $\lambda$ in $\mathbb{C}$ of modulus $1$ such that
$h_\phi(z)^2\phi(z)=\lambda z\phi'(z)^2$ for all $z$ in $\Delta$. We observe 
that the right-hand side has a zero of odd degree at $0$. Hence $h_\phi(z)^2
\phi(z)$ has must have a zero of odd degree at $0$. In particular, we have that
$\phi(0)=0$. 

Let us write $\phi$ as $\phi(z)=z\psi^2(z)$ where $\psi\colon 
\Delta\to\Delta$ is analytic. Note that since $h_\phi(z)^2\phi(z)=\lambda 
z\phi'(z)^2$ for all $z$ in $\Delta$ and $h_\phi$ is non-zero on ${ \mathcal 
B}_v^T
(\Delta)$ we have that each zero of $\psi$ in ${ \mathcal B}_v^T(\Delta)$ is of
order $1$. Moreover, for $z$ in $\Delta$ we have that 
$$
h_\phi(z)^2=\lambda\frac{z}{\phi
(z)}\phi'(z)^2=\lambda(\psi(z)+2z\psi'(z))^2.
$$
 Using \cite[Theorem~17.9]{Rudin} write 
$\phi$ as
$\phi(z)=\lambda z^kB(z)^2g(z)^2$ where $k\in \mathbb{N}$, $B(z)$ is the 
Blaschke product formed by the roots of $\psi$ and $g$ is a non-zero bounded 
holomorphic function on $\Delta$ with $\|g\|_\infty=\|\psi\|_\infty=1$. Suppose 
that $a\not=0$ is a zero of $\phi$. Write $B(z)$ as 
$$
B(z)=\frac{|a|}{a}\frac{z-a}{1-\bar{a}z}B_a(z).
$$
Then we have that 
$$
\frac{1}{1-|a|}=|h_\phi(a)|=2|a\psi'(a)|=2|a|\left(\frac{|a|^{\frac{k-1}{2}}}
{1-|a|^2}\right)|B_a(a)g(a)|
$$
giving that 
$$|B_a(a)g(a)|=\frac{1+|a|}{2|a|^{\frac{k+1}2}}$$
which is impossible as $|a|<1$ and $\|B_ag\|_\infty\leqslant 1$. Hence we have 
that $0$ is the only root of $\phi$ in ${\mathcal B}_v^T(\Delta)$.  
Repeating the argument of the case where
$1<\beta<2$ we get that $\phi$ must be an automorphism of the disc. This 
completes the proof.

\begin{theorem}\label{beta2}
Let $\beta\geqslant 1$ and $v\colon\Delta\to\Delta$ be given by 
$\displaystyle{v(z)=e^{\frac{-1}{1-|z|^\beta}}}$. Then every isometry $T\colon
{\mathcal H}_{v_o}(\Delta)\to 
{ \mathcal H}_{v_o}(\Delta)$ has the form $T(f)(z)=\alpha f(e^{i\theta}z)$, $z\in 
\Delta$, for some complex number $\alpha$ with $|\alpha|=1$, and some $\theta\in
\mathbb R$. In particular, $T$ is surjective.
\end{theorem}

{\sc Proof:} From  Theorem~\ref{stone} we know that there is an 
analytic function $\phi\colon
\Delta\to\Delta$ and $h_\phi$ in ${ \mathcal H}_{v_o}(\Delta)$ such that $(Tf)(z)=
h_\phi(z)f\circ\phi(z)$ for all $z$ in $\Delta$. Moreover for each $z$ in 
${ \mathcal B}_v^T(\Delta)$ we have that
$$|h_\phi(z)|=\frac{\exp^{\frac{-1}{1-|\phi(z)|^\beta}}}{\exp^{\frac{-1}{1-|z|^\beta}}}.$$
Taking Laplacians of $\log |h_\phi|$ we get that
$${(|\phi(z)|^{2(\beta-1)}+|\phi(z)|^{\beta-2})|\phi'(z)|^2\over (1-|\phi(z)
|^\beta)^3}={(|z|^{2(\beta-1)}+|z|^{\beta-2})\over (1-|z|^{\beta})^3}.$$
By Theorem \ref{interior} for each $r>0$, 
$\overline{{ \mathcal B}_v^T(\Delta)\cap \phi^{-1}(B(0,r))}$ has non-empty 
interior.
This means that for each $n$ in $\mathbb{N}$ we can choose $a_n$ in the 
interior of $\overline{{ \mathcal B}_v^T(\Delta)\cap\phi^{-1}
(B(0,\frac{1}{n}))}$.
As in the case with $v(z)=(1-|z|^\beta)$ we may suppose that $(a_n)_n$ 
converges
to a point $a$ of $\Delta$. Moreover, for each $n$ in $\mathbb{N}$ we have that
$${(|\phi(a_n)|^{2(\beta-1)}+|\phi(a_n)|^{\beta-2})|\phi'(a_n)|^2\over (1-|\phi(a_n)
|^\beta)^3}={(|a_n|^{2(\beta-1)}+|a_n|^{\beta-2})\over (1-|a_n|^{\beta})^3}.$$

We will again distinguish between four different values for $\beta$. First we 
consider the case when $\beta>2$. In this case we see that the sequence 
$(a_n)_n$ chosen above is a null sequence and by continuity $\phi(0)=0$. For 
$z$ in the interior of $\overline{{ \mathcal B}_v^T(\Delta)}$ we have that 
$$|\phi'(z)|^2=\frac{|z|^{2(\beta-1)}+|z|^{\beta-2}}{|\phi(z)|^{2(\beta-1)}
+|\phi(z)|^{\beta-2}}\frac{(1-|\phi(z)|^\beta)^3}{(1-|z|^\beta)^3}.$$
The Schwarz lemma implies that the right-hand side is greater than or equal to
$1$. Letting $z$ tend to $0$ we get that $|\phi'(0)|=1$ and hence we have that
$\phi$ is an automorphism. 

When $\beta=2$ then for any $z$ in the interior of $\overline{{ \mathcal B}_v^T(
\Delta)}$ we have that 
$${(|\phi(z)|^2+1)|\phi'(z)|^2\over (1-|\phi(z)|^2)^3}={|z|^{2}+1
\over (1-|z|^2)^3}.$$
It follows from the Schwarz-Pick Lemma that 
$$\frac{1+|z|^2}{1-|z|^2}\leqslant \frac{1+|\phi(z)|^2}{1-|\phi(z)|^2}$$
which gives that $|z|\leqslant |\phi(z)|$ for any $z$ in the interior of 
$\overline{{ \mathcal B}_v^T(\Delta)}$. In particular, we get 
that $a=0$. The Schwarz Lemma now implies that $\phi$ is an automorphism and
hence $T$ is surjective.

Let us now consider the case when $1\leqslant\beta<2$. Suppose that 
$a\not=0$. Then let $k$ be the degree of the zero of $\phi$ at $a$. As 
$${(|\phi(z)|^{2(\beta-1)}+|\phi(z)|^{\beta-2})|\phi'(z)|^2\over (1-|\phi(z)
|^\beta)^3}={(|z|^{2(\beta-1)}+|z|^{\beta-2})\over (1-|z|^{\beta})^3}$$
we see that $\displaystyle{\frac{|\phi'(z)|^2}{|\phi(z)|^{2-\beta}}}$ must gave a 
finite non-zero limit at $a$. Hence $2(k-1)=(2-\beta)k$ or $k\beta=2$. However 
if $1<\beta<2$ this is impossible and thus $a=0$ is the only 
zero of $\phi$ in ${ \mathcal B}_v^T(\Delta)$.

If $\beta=1$ then $\phi$ has
a double zero at $a$. Writing $\phi$ as $\phi(z)=(z-a)^2\psi(z)$ we see that 
$$\frac{\left(1+\frac{1}{|(z-a)^2\psi(z)|}\right)}{(1-|\phi(z)|)^3}\left|
(z-a)^2\psi'(z)+2(z-a)\psi(z)\right|^2=\frac{\left(1+\frac{1}{|z|}\right)}
{(1-|z|)^3}.$$
Setting $z=a$ we get that
$$4|\psi(a)|=\frac{\left( 1+\frac{1}{|a|}\right)}{(1-|a|)^3}.$$  
We have that 
$$|\psi(a)|\leqslant \max_{|z|=1}\frac{|\phi(z)|}{|z-a|^2}\leqslant 
\frac{1}{(1-|a|)^2}.$$
Therefore 
$$\frac{\left( 1+\frac{1}{|a|}\right)}{(1-|a|)^3}\leqslant\frac{4}{(1-|a|)^2}$$
or 
$$\frac{\left(1+\frac{1}{|a|}\right)}{(1-|a|)}\leqslant 4.$$
However as the minimum value of the function $f(r)=\displaystyle{\frac{
\left(1+\frac
{1}{r}\right)}{(1-r)}}$ over the internal $(0,1)$ is $\displaystyle{\frac{1+
\frac{1}{-1+\sqrt{2}}}{2-\sqrt{2}}}$ which is equal to $5.8284...$ we have a
contradiction and thus $a=0$ is the only zero of $\phi$ in ${ \mathcal B}_v^T
(\Delta)$. 

Returning to the case where $1\leqslant\beta<2$, for each $n$ in 
$\mathbb{N}$ we have 
$${(|\phi(a_n)|^{2(\beta-1)}+|\phi(a_n)|^{\beta-2})|\phi'(a_n)|^2\over (1-|\phi(a_n)
|^\beta)^3}={(|a_n|^{2(\beta-1)}+|a_n|^{\beta-2})\over (1-|a_n|^{\beta})^3}.$$
Multiplying by $|a_n|^{2-\beta}$ and letting $n$ tend to infinity to get that
$|\phi'(0)|=1$ which means that $\phi$ is an automorphism of the disc and hence
$T$ is surjective.

\begin{theorem}\label{beta3}
Let $\beta<0$ and $v\colon\Delta\to\Delta$ be given by $v(z)=(1-\log(1-|z|)
)^\beta$. Then every isometry $T\colon { \mathcal H}_{v_o}(\Delta)\to 
{ \mathcal H}_{v_o}(\Delta)$ has the form $T(f)(z)=\alpha f(e^{i\theta}z)$, $z\in 
\Delta$, for some complex number $\alpha$ with $|\alpha|=1$, and some $\theta
\in\mathbb R$. In particular $T$  is surjective.
\end{theorem}

{\sc Proof:} We have that $\Delta \log v(z)$ is given by 
\begin{align*}
\Delta(\log(v(|z|))=&\beta \left( {1\over |z|(1-|z|)^2(1-\log(1-|z|))}-
{1\over (1-|z|)^2(1-\log(1-|z|))^2}\right)\\
\end{align*}
The result now follows in the same way as with the weight 
$\displaystyle{e^{\frac{-1}{(1-|z|^\beta)}}}$.

We have seen that for the weights $v(z)=1-|z|$ and $w(z)=1-|z|^2$ on the unit 
disc $\Delta$ the mapping $T(f)(z)=f(z^2)$ is an example of a 
non-surjective isometry from ${ \mathcal H}_{v_o}(\Delta)$ into ${ \mathcal 
H}_{w_o}(
\Delta)$. Let us now observe that, up to composition with an automorphism of
the disc, all isometries between these spaces are of this form.

\begin{theorem}
Let $v(z)=1-|z|$ and $w(z)=1-|z|^2$ on $\Delta$. 
Then every isometry $T$ from ${ \mathcal H}_{v_o}(\Delta)$ into ${\mathcal H
}_{w_o}
(\Delta)$ has the form $T(f)(z)=\psi'(z)f(\psi(z)^2)$ for some automorphism 
$\psi$ of $\Delta$. 
\end{theorem}

{\sc Proof:} Let $T$ be an isometry from ${ \mathcal H}_{v_o}(\Delta)$ 
into 
${ \mathcal H}_{w_o}(\Delta)$. Then, by Theorem~\ref{stone}, 
there is an analytic function $\phi\colon\Delta\to
\Delta$ and $h_\phi$ in ${ \mathcal H}_{w_o}(\Delta)$ such that $T(f)(z)=h_\phi(z)f
\circ\phi(z)$ for all $z$ in $\Delta$. Moreover for each $z$ in $\overline{
{ \mathcal B}_w^T(\Delta)}$ we have that
$$|h_\phi(z)|=\frac{1-|\phi(z)|}{1-|z|^{2}}.$$
Since $h_\phi$ is analytic we have that $\log |h_\phi(z)|$ is harmonic on 
$\Delta\setminus h_\phi^{-1}(0)$. Hence for $z$ in the interior of $\overline{
{ \mathcal B}_v^T(\Delta)}$ we have that
$$\Delta\log(1-|\phi(z)|)=\Delta\log(1-|z|^{2}).$$ 
This gives us that 
$$\frac{|\phi(z)|^{-1}|\phi'(z)|^2}{(1-|\phi(z)|)^2}=\frac{4}{(1-|z|^{2})^2}.$$

We rewrite the above equation as 
$$\left|\frac{(1-|\phi(z)|)^2}{(1-|z|^2)^2}\right|=\frac{|\phi'(z)^2|}
{4|\phi(z)|}$$
or as 
$$|h_\phi(z)^2|=\frac{|\phi'(z)^2|}{4|\phi(z)|}$$
for all $z$ in the interior of $\overline{{ \mathcal B}_v^T(\Delta)}$. As in 
Theorem~\ref{beta1} we get that there is $\lambda$ in $\mathbb{C}$ with 
$|\lambda|=1$ so that $h_\phi(z)^2=\lambda \frac{\phi'(z)^2}{4\phi(z)}$ for all 
$z$ in $\Delta$. In
particular, we observe that there is an analytic function $\psi\colon \Delta
\to \Delta$ such that $\phi(z)=\psi(z)^2$ for all $z$ in $\Delta$. Taking $z$ in
the interior of $\overline{{ \mathcal B}_w^T(\Delta)}$ and replacing $\phi$ 
with 
$\psi^2$ we see that 
$$|\psi'(z)|^2=\frac{(1-|\psi(z)|^2)^2}{(1-|z|^2)^2}$$
for all $z$ in $\Delta$. 
The Schwarz-Pick Lemma now implies that $\psi$ is an automorphism of $\Delta$ 
and the result follows.

\section{Isometries of ${ \mathcal H}_v(U)$}

Let us start with an 
example of a non-surjective isometry of ${\mathcal H}_v(U)$.

\begin{example}\label{exbloch}
The following example is due to Bonet, Lindstr\"om and Wolf 
\cite{BLW} which in turn
is motivated by an example of a non-surjective isometry of the little Bloch 
space given by Martin and Vukoti\'c \cite{MV}. 
Given a thin interpolating sequence $(a_n)_n$ and a non-negative integer $m$ 
we form the corresponding Blaschke product
$$B(z)=z^m\prod_{n=1}^\infty \frac{\overline{a_n}}{|a_n|}\frac{z-a_n}{1-z
\overline{a_n}}.$$
An interpolating sequence 
$(a_n)_n$ in $\Delta$ with $a_n\not= 0$ for all $n$ is said to be thin if  
$$\lim_{n\to\infty}\prod_{k\not= n}\left|\frac{a_k-a_n}{1-a_n\overline{a_k}}\right|
=\lim_{n\to \infty} (1-|a_n|^2)|B'(a_n)|=1.$$
Consider the weight $v(z)=1-|z|^2$ on the unit disc $\Delta$. Then
$T_B$ given by
$$T_B(f)(z)=B'(z)f\circ B(z)$$
is a non-surjective isometry from ${ \mathcal H}_v(\Delta)$ into ${ \mathcal H}_v
(\Delta)$. See \cite{BLW} and \cite{MV} for the details.
\end{example}

So, for weights such as $v(z)=1-|z|^2$ there are non-surjective isometries
from ${ \mathcal H}_v(\Delta)$ into ${ \mathcal H}_v(\Delta)$. The following result shows
that such isometries can be characterised by how they map ${ \mathcal H}_{v_o}(
\Delta)$.

\begin{theorem}
Consider the weights $v_1(z)=1-|z|^\beta$, $\displaystyle{v_2(z)=e^{\frac{-1}{1-
|z|^\beta}}}$ $\beta\geqslant 1$ and the weight
$v_3(z)=(1-\log(1-|z|))^\beta$, $\beta<0$, on the open unit disc $\Delta$. Let
$T\colon { \mathcal H}_{v_i}(\Delta)\to{ \mathcal H}_{v_i}(\Delta)$ be an 
isometry 
($i=1,2$ or $3$). Then $T$ is surjective if and only if $T\left({ \mathcal 
H}_{(v_i)_o}(\Delta)\right)\subseteq { \mathcal H}_{(v_i)_o}(\Delta)$.
\end{theorem}

{\sc Proof:} Let us just write $v$ for $v_1$, $v_2$ or $v_3$. Suppose 
that $T\colon { \mathcal H}_{v}(\Delta)\to{ \mathcal 
H}_{v}(\Delta)$ is surjective. Since ${ \mathcal H}_{v_o}(\Delta)$ is an M-ideal
in ${ \mathcal H}_v(\Delta)$ \cite[Theorem~4.2]{HL} implies that $T$ is the 
bitranspose of the isometric isomorphism $T|_{{ \mathcal H}_{v_o}(\Delta)} \colon 
{ \mathcal H}_{v_o}(\Delta)\to { \mathcal H}_{v_o}(\Delta)$ and so maps 
${ \mathcal H}_{v_o}(
\Delta)$ into ${ \mathcal H}_{v_o}(\Delta)$. Conversely, if 
$T|_{{ \mathcal H}_{v_o}(\Delta)} \colon 
{ \mathcal H}_{v_o}(\Delta)\to { \mathcal H}_{v_o}(\Delta)$ then, by 
Theorem~\ref{beta1}, Theorem~\ref{beta2} and Theorem~\ref{beta3}, it 
is surjective and therefore $T=(T|_{{ \mathcal H}_{v_o}(\Delta)})''$ is a surjective 
isometry
from ${ \mathcal H}_v(\Delta)$ onto ${ \mathcal H}_v(\Delta)$.

If we consider the example of Bonet, Lindstr\"om and Wolf of a non-surjective 
isometry $T_B$ for the weight $v(z)=1-|z|^2$ (see Example~\ref{exbloch}) we have
that $T_B(1)(z)=B'(z)$ where $B$ is a thin Blaschke product. Then $\lim_{n\to
\infty}(1-|a_n|^2)|B'(a_n)|=1$ and therefore $T_B(1)$ belongs to 
${ \mathcal H}_v(\Delta)$ but not to ${ \mathcal H}_{v_o}(\Delta)$.

\section{Isometries of the Bloch Space}

In \cite{MV} Mart{\'\i}n and Vukoti{\' c} use the hyperbolic derivative and 
cluster sets to characterise the isometric composition operators between the 
Bloch space of all holomorphic functions $f\colon \Delta\to \mathbb{C}$ such 
that $\|f\|=|f(0)|+\sup_{z\in\Delta}(1-|z|^2)|f'(z)|<\infty$.

The normalised Bloch, ${ \mathcal B}$, is defined as the space of holomorphic 
functions $f\colon \Delta\to \mathbb{C}$ such that $f(0)=0$ and 
$\|f\|_{ \mathcal B}:=
\sup_{z\in\Delta}(1-|z|^2)|f'(z)|<\infty$. The little Bloch space is the set of 
all $f$ in ${ \mathcal B}$ such that $\lim_{|z|\to 1}(1-|z|^2)|f'(z)|
=0$ and is denoted by ${\mathcal B}_o$. Setting $v(z)=1-|z|^2$ we see that the 
mapping $D\colon f\mapsto f'$ is an isometric isomorphism of ${ \mathcal B}$ 
onto ${ \mathcal 
H}_v(\Delta)$ which maps ${ \mathcal B}_o$ onto ${ \mathcal H}_{v_o}(\Delta)$.  
 Using this
identification Cima and Wogen, \cite{CW}, showed that all isometries of the 
little Bloch space, ${ \mathcal B}_o$, are surjective (see also \cite{FJ1}).
Their proof (and so the one in \cite{FJ1}) unfortunately seems to be 
incomplete. Both proofs show that there is a 
subset $\Sigma({\mathcal R}_0)$ of $\Delta$ and functions $\tau\colon \Sigma
({\mathcal R}_0)\to \Delta$, $\alpha\colon \Sigma({\mathcal R}_0)\to \Gamma$ 
such that $T^*(\delta_z)=\alpha(z)\delta_{\tau(z)}$. The space $\Sigma({\mathcal 
R}_0)$ may be regarded as 
corresponding to our ${\mathcal B}_w^T(\Delta)$ while $\tau$ corresponds to our
$\phi_1$. It is then shown 
that there is a holomorphic function $G_0$ on the unit disc so that 
$|G_0(z)|=\frac{1-|\tau(z)|^2}
{1-|z|^2}$
for $z$ in $\Sigma({\mathcal R}_0)$. The function $\tau$ is extended to a 
holomorphic function of the disc into the disc. At the end of the theorem, 
Lemma~1 of \cite{CW} is applied to $G_0$ and then it is concluded that $\tau$ 
is an automorphism of 
the disc. However in order to apply Lemma~1 as stated the equality 
$|G_0(z)|=\frac{1-|\tau(z)|^2}
{1-|z|^2}$
should hold on $\Delta$. (A look at the proof however shows that an open subset 
of $\Delta$ would suffice.) As far as we can see this equality only occurs on
$\Sigma({\mathcal R}_0)$. In addition, in order to classify the isometries 
of the Little Bloch space with 
the norm $\||f|\|=\sup_{z\in\Delta}(1-|z|)|f'(z)|$ the application of
Lemma~1' requires that there is a function $\tau\colon\Delta\to\Delta$ which
satisfies $|f(z)|=\frac{1-|\tau(z)|}{1-|z|}$
for some analytic function $f$ on $\Delta$. However, we can only see that this
equality will hold on a distinguished subset of $\Delta$. 

However, using Theorem~\ref{interior} and 
Theorem~\ref{beta1} with $\beta=2$, we are able to recover \cite[Theorem 1]
{CW} and show that all isometries of the little Bloch space are indeed
surjective. Theorem~\ref{beta1} also shows that each isometry of the 
(normalised) little
Bloch space with the norm $\|f\|_\alpha=\sup_{z\in\Delta}(1-|z|^\alpha)
|f'(z)|$, $1\leqslant\alpha<\infty$, is surjective.

Setting $Z$ equal to the set $\{(1-|z|^2)f':f\in { \mathcal 
B}\}$ we see that each isometry $T\colon { \mathcal B}\to { \mathcal B}$ 
induces an 
isometry $\tilde T\colon Z\to Z$ such that the following diagram commutes
$$
\xymatrix{{ \mathcal B}\ar[r]^T \ar[d]_D& { \mathcal B}\ar[d]_D & \\ 
Z\ar[r]^{{\tilde T}}\ar[r]& Z.\\}
$$ 
As in the case with ${ \mathcal H}_v(U)$ the extreme points of the closed unit 
ball
of $({\tilde T}(Z))'$ are of the form $\delta_z$ for $z$ in $\beta \Delta$, 
the Stone--{\v C}ech compactification of the unit disc. We set ${ \mathcal B}^T$
 equal to the set of all $z$ in $\beta\Delta$ such that $\delta_z$ 
is an extreme point of the closed unit ball of $({\tilde T}(Z))'$. As with 
our previous Banach--Stone Theorems we have that $\tilde T$ induces a 
function $\phi$ from ${ \mathcal B}^T$ onto $\beta\Delta$ and such that 
${\tilde T}f(z)=\lambda f\circ\phi(z)$ for all $f$ in $Z$ and all $z$ in 
$\Delta\cap {\mathcal B}^T$ some $\lambda$ in $\mathbb{C}$ with $|\lambda|=1$.
The isometry $f\mapsto f'$ and our criteria for surjectivity of the
weighted space ${ \mathcal H}_v(\Delta)$ gives the following result.

\begin{theorem}
Let $T\colon { \mathcal B}\to { \mathcal B}$ be an isometry. Then the 
following are equivalent
\begin{enumerate}
\item[(a)] $T$ is surjective,
\item[(b)] $T({ \mathcal B}_o)\subseteq { \mathcal B}_o$,
\item[(c)] $\Delta\cap \phi_1^{-1}(\Delta)$ has non-empty interior.
\end{enumerate}
\end{theorem}
   
{\bf Acknowledgements:}\label{ackref}
The authors wish to thank Stephen Gardiner for his suggestions concerning
Theorem~\ref{interior}, Richard Smith for his reference to the Invariance
of Domains and Joseph Cima for his correspondence regarding 
\cite{CW}.

Christopher Boyd \\
School of Mathematical Sciences,\\
University College Dublin,\\
Belfield,\\
Dublin 4,\\ Ireland.
email:{Christopher.Boyd@ucd.ie}

\bigskip
Pilar Rueda\\
Departamento de An\'alisis Matem\'atico,\\
Facultad de Matem\'aticas,\\
Universidad de Valencia,\\
46100 Burjasot,\\
Valencia,\\
Spain.
 email:{Pilar.Rueda@uv.es}


\begin{thebibliography}{99}

\bibitem{AF} {\sc J. Araujo \& J. Font.} Linear isometries between
subspaces of continuous functions, {\em Trans. Amer. Math. Soc.,}
{\bf 349} (1), (1997), 413--428.

\bibitem{BiSu} {\sc K.D. Bierstedt \& W.H. Summers.} Biduals of weighted Banach 
spaces  of analytic functions, {\em J. Austral. Math. Soc.,} {\bf 54} (1993), 
70--79.

\bibitem{BLW} {\sc J. Bonet, M. Lindstr\"om \& E. Wolf.}
Isometric weighted composition operators on weighted Banach spaces of type 
${ \mathcal H}^\infty$, {\em Proc. Amer. Math. Soc.,} {\bf 136} (2008), no. 12, 
4267--4273.
 
\bibitem{BR1} {\sc C. Boyd \& P. Rueda.} The $v$-boundary of
weighted spaces of holomorphic functions, {\em Ann. Acad. Sci. Fenn.
Math.,} {\bf 30} (2005), 337--352.

\bibitem{BR2} {\sc C. Boyd \& P. Rueda.} Complete weights and
$v$-peak points  of spaces of weighted holomorphic functions,
{\em Israel Journal of Mathematics,} {\bf 155} (2006), 57-80.

\bibitem{BR3} {\sc C. Boyd \& P. Rueda.} Bergman and Reinhardt
weighted spaces of holomorphic functions, {\em Illinois J. Math.,}
{\bf 49} (1), (2005), 217--236.

\bibitem{BR4} {\sc C. Boyd \& P. Rueda.} Isometries between spaces of
weighted holomorphic functions, {\em Studia Math.,} {\bf 190} (3),
(2009), 203--231.

\bibitem{BR5} {\sc C. Boyd \& P. Rueda.} Isometries of weighted spaces
of holomorphic functions on unbounded domain, {\em Proc. Roy. Soc. Edinburgh
Sect. A,} {\bf 139} (2009), 253--271.

\bibitem{BR6} {\sc C. Boyd \& P.Rueda.} The biduality problem and M-ideals in 
weighted spaces of holomorphic functions. {\em J. Convex Anal.,} {\bf 18} (4), 
(2011),  1065--1074.

\bibitem{CW} {\sc J.A. Cima \& W.R. Wogen.} On isometries of the Bloch 
space, {\em Illinois J. Math.,} {\bf 24} (2), (1980), 313--316.
 
\bibitem{DuSc} {\sc N. Dunford \& J.T. Schwartz,} Linear Operators. Part I: 
General Theory. {\em John Willey \& Sons,} 1957. 

\bibitem{FJ1} {\sc J. Fleming \& J.E. Jamison,.} Isometries on Banach 
spaces: function spaces. Chapman \& Hall/CRC Monographs and Surveys in Pure and 
Applied Mathematics, 129. {\em Chapman \& Hall/CRC, Boca Raton, FL,} 2003.

\bibitem{FJ2} {\sc J. Fleming \& J.E. Jamison.} Isometries on Banach 
spaces. Vol. 2. Vector-valued function spaces. Chapman \& Hall/CRC Monographs 
and Surveys in Pure and Applied Mathematics, 138. {\em Chapman \& Hall/CRC, 
Boca Raton, FL,} 2008.  

\bibitem{HL} {\sc P. Harmand \& A. Lima.} Banach spaces which are
M-deals in their biduals, {\em Trans. Amer. Math. Soc.,} {\bf 283} (1984), 
253--264.

\bibitem{HWW} {\sc P. Harmand, D. Werner \& W. Werner.} M-Ideals in Banach
spaces and Banach Algebras, {\em Lecture Notes in Mathematics,} vol. 1547 
Springer,(1993).

\bibitem{Lu1} {\sc W, Lusky.} On weighted spaces of harmonic and holomorphic  
functions, {\em J. London Math. Soc.,} {\bf 51} (1995), 309--320.

\bibitem{Lu2} {\sc W. Lusky.} On the isomorphism classes of weighted spaces of 
harmonic and holomorphic functions. {\em Studia Math.,} {\bf 175} (1), (2006), 
19--45.

\bibitem{MV} {\sc M.J. Martin \& D. Vukotic.} Isometries of the Bloch 
space among the composition operators, {\sl Bull. London Math. Soc.,} {\bf 39} 
(2007), 151-155.

\bibitem{Mu} {\sc J. Mujica.}  Complex analysis in Banach spaces. 
Holomorphic functions and domains of holomorphy in finite and infinite 
dimensions, North-Holland Mathematics Studies, vol. 120  {\em North-Holland 
Publishing Co., Amsterdam,} 1986.

\bibitem{MP} {\sc A. Mukherjea \& K. Pothoven.} Real and Functional 
Analysis, Part A, Real Analysis, Second Edition, {\em Mathematical Concepts 
and Methods in Science and Engineering,} vol. 27, Plenum Press, New York and 
London, (1984).
   
\bibitem{Rudin} {\sc W. Rudin.} Real and Complex Analysis, {\em 
McGraw-Hill,} (1987).
\end{thebibliography}
\end{document}